\documentclass[11pt]{article}
\usepackage{amssymb,amsfonts,amsmath}%
\usepackage{latexsym}
\usepackage{graphicx,epsfig}
\usepackage{color}

\newtheorem{thm}{Theorem}
\newtheorem{lem}[thm]{Lemma}
\newtheorem{pro}[thm]{Proposition}
\newtheorem{cor}[thm]{Corollary}

\newcommand{\e}{\emph}
\newcommand{\al}{\alpha}
\newcommand{\tho}{\widetilde{H}_1}
\newcommand{\tht}{\widetilde{H}_2}
\newcommand{\thi}{\widetilde{H}_i}

\newcommand{\hyp}{$k$-unifrom hypergraph }
\newcommand{\hyps}{$k$-uniform hypergraphs }
\newcommand{\hypse}{$k$-unifrom hypergraphs. }

\newcommand {\cbdo}{\hfill$\Box$}

\begin{document}
\title{\bf A note on packing of uniform hypergraphs}
\author{Jerzy Konarski, Mariusz Wo\'zniak and Andrzej \.{Z}ak\thanks{The authors were partially supported
by the Polish Ministry of Science and Higher Education (11.11.420.004).}\\
\small{AGH University of
Science and Technology, 
Al. Mickiewicza 30, 30-059 Krak\'ow,
Poland}} \maketitle
\begin{abstract}
A \e{packing} of two $k$-uniform hypergraphs $H_1$ and $H_2$ is a set $\{H_1', H_2'\}$ of edge-disjoint sub-hypergraphs of the complete $k$-uniform hypergraph $K_n^{(k)}$ such that $H_1'\cong H_1$ and $H_2'\cong H_2$. Whilst the problem of packing of graphs (i.e. 2-uniform hypergraphs) has been studied extensively since seventies with many sharp results, much less is known about packing of general hypergraphs. In this paper we attempt to find the minimum possible sum of sizes $m(n,k)$ of two $k$-uniform, $n$-vertex hypergaphs which do not pack. We also prove a sufficient condition on the product of maximum degrees, which guarantees the packing.
\end{abstract}

\section{Introduction}
By a \e{hypergraph} $H$ we mean a pair $(V(H), E(H))$ where $V(H)$ is a finite set (elements of $V(H)$ are called \e{vertices}) and $E(H)$ is a family of subsets of $V(H)$ (members of $E(H)$ are called \e{edges}). We use the term {\it $k$-uniform hypergraph} to refer to hypergraphs such that each edge consists of exactly $k$ vertices.

Two \hyps $H_1$ and $H_2$ \e{pack} (into a complete $k$-uniform hypergraph $K_n^{(k)}$, where $n=\max\{|V(H_1)|,|V(H_2)|\}$) if there is a pair of edge-disjoint subhypergraphs $\{H_1',H_2'\}$ of $K_n^{(k)}$ 
such that $H_1\cong H_1'$ and $H_2\cong H_2'$. The problem of packing of graphs (i.e. 2-uniform hypergraphs) has been considered by many authors since in 1978 Bollob\'as and Eldridge \cite{BE} and Sauer and Spencer \cite{SS} proved first important results.
In particular Sauer and Spencer \cite{SS} showed the following theorems. 
\begin{thm}[\cite{SS}]\label{SS4} Let $G_1$ and $G_2$ be two graphs on $n$ vertices. If
\[|E(G_1)|\cdot|E(G_2)|<\binom{n}{2}\]
then $G_1$ and $G_2$ pack.
\end{thm}
\begin{thm}[\cite{SS}]\label{SS2} Let $G_1$ and $G_2$ be two graphs on $n$ vertices. If
\[2\Delta(G_1)\Delta(G_2)<n,\]
then $G_1$ and $G_2$ pack.
\end{thm}
\begin{thm}[\cite{SS}]\label{SS1} Let $G_1$ and $G_2$ be two graphs on $n$ vertices. If
\[|E(G_1)|+|E(G_2)|\leq \left\lceil\frac{3}{2}n\right\rceil-2,\]
then $G_1$ and $G_2$ pack. 
\end{thm}

Our purpose is to find some analogues of these theorems in the case 
of $k$ uniform hypergraphs, $k\geq 3$. 
Theorem \ref{SS4} has been generalized by Naroski \cite{N}.   
\begin{thm}[\cite{N}]\label{nar} Let $H_1$, $H_2$ be two \hyps on $n$ vertices. If
\[|E(H_1)|\cdot|E(H_2)|<\binom{n}{k}\] then $H_1$ and $H_2$ pack.
\end{thm}
(In fact Theorem \ref{nar} is a special case of a more general result from \cite{N}).

A generalization of Theorem \ref{SS2} has been obtained by R\"odl, Ruci\'nski and Taraz \cite{RRT}. 
Let $d_{l}(U)=|\{e\in E(H)\colon U\subset e\}|$ be a degree of an $l$-element subset $U$ of vertices in a hypergraph $H$ and $\Delta_{l}(H)=\max\{d_l(U)\colon U\in \binom{V(H)}{l}\}$. Let $\Delta_1=\Delta$ as in the simple graphs. 
\begin{pro}[\cite{RRT}]\label{lemma}
Let $H_1$ and $H_2$ be \hyps on $n$ vertices. If
\[\Delta(H_1)\Delta_{k-1}(H_2)+\Delta(H_2)\Delta_{k-1}(H_1)<n-k+2.\]
then $H_1$ and $H_2$ pack.
\end{pro}
R\"odl, Ruci\'nski and Taraz  used this proposition in order to obtain 
a far reaching improvement in the case when one hypergraph has bounded maximum degree. 
We further generalize Theorem \ref{SS2} in the following way.
\begin{thm} \label{beta}
Let $H_1$ and $H_2$ be \hyps on $n$ vertices. If there exists $\beta$, $0<\beta<k$ such that
\begin{align}\label{beta}\Delta_{\beta}(H_1)\Delta_{k-\beta}(H_2)+\Delta_{k-\beta}(H_1)\Delta_{\beta}(H_2)<\binom{n}{\beta}-\binom{k}{\beta}+2\end{align}
then $H_1$ and $H_2$ pack.
\end{thm}

Let us note, that the bound in Proposition \ref{lemma} is far from being convenient as the term $\Delta(H_i)$ for $i=1,2$ may be as large as $\binom{n}{k-1}$ (as the authors of \cite{RRT} admit). This drawback does not exist in the case $\beta=\frac{k}{2}$ (for even $k$) in Theorem \ref{beta}.

Finally, we attempt to find a (tight) bound on the sum of sizes of two \hyps that guarantees 
the existence of a packing. To this end we define $m(n,k)$ to be the minimum possible number $m$ such that 
there exist two $k$-uniform hypergaphs $H_1$ and $H_2$ which do not pack 
with $|E(H_1)|+|E(H_2)|=m$ and $\max\{|V(H_1)|,|V(H_2)|\}=n$. Thus, Theorem \ref{SS1} yields that 
\[m(n,2)=\left\lceil\frac{3}{2}n\right\rceil-1\]
This bound has been considerably weakened by Bollob\'as and Eldridge \cite{BE} in the case when neither graph has a total vertex, 
which in turn has been extended to general hypergraphs without edges of size 0,1,$n-1$, $n$ by  
Kostochka, Stocker and Hamburger in \cite{KSH} (weaker versions appeared earlier in \cite{N,PW1,PW2}). One can expect further improvements for 
$k$-uniform hypergraphs for $k\geq 3$. 

An immediate corollary of Theorem \ref{SS4} is the following. 
\begin{cor}
Let $H_1$, $H_2$ be \hyps on $n$ vertices. If 
\[|E(H_1)|+|E(H_2)|<2\sqrt{\binom{n}{k}},\] then $H_1$ and $H_2$ pack.
\end{cor}

This implies that 
\begin{align}\label{order}
m(n,k)=\Omega(n^{k/2}).
\end{align}
Our second result shows that the order of magnitude in (\ref{order}) is correct if 
$k$ is even and assymptotically correct if $k$ is odd.

\begin{thm} \label{final}
Let $k=2\alpha$ be even. Let $n$ be such that $\binom{k-i}{k/2-i}$ divides $\binom{n-i}{k/2-i}$ for all $i$, $0\leq i \leq k/2-1$. Then
\[m(n,k)\leq \binom{n-\al}{\al}+\frac{\binom{n}{\al}}{\binom{2\al}{\al}}.\]
Furthermore, for odd $k$ and $n$ such that $k-i$ divides $\binom{n-i}{k-1-i}$ for all $i$,  $0\leq i \leq k-2$.
\[m(n,k)\leq cn^{(k^2-k-1)/(2k-3)},\]
where $c$ is some constant depending only on $k$.
\end{thm}

\section{Proofs}
A $t-(n,k,\lambda)$-design on a set $X$ of size $n$ is a collection $T$ of $k$-element subsets of $X$ such that every $t$ elements of $X$ are contained in exactly $\lambda$ sets in $T$. It is easy to observe that a necessary condition for the existence of a $2$-$(n,3,1)$ design (called Steiner Triple System) is that $n$ must be congruent to 1 or 3 mod 6. In 1846, Kirkman showed that this necessary condition is also sufficient. In 1853, Steiner posed the natural generalisation of the question: given $q$ and $r$, for which $n$ is it possible to choose a collection $Q$ of $q$-element subsets of an $n$-element set $X$ such that any $r$ elements of $X$ are contained in exactly one of the sets in $Q$? There are some natural necessary divisibility conditions generalising the necessary conditions for Steiner Triple Systems. The Existence Conjecture states that for all but finitely many $n$ these divisibility conditions are
also sufficient for the existence of general Steiner systems (and more generally
designs). Recently this conjecture was proved by Keevash \cite{K}.

\begin{thm}[\cite{K}]\label{Keevash}
A $t-(n,k,\lambda)$-design on a set $X$ exists if and only if for every $0\leq i \leq t-1$
\[\binom{k-i}{t-i}\quad\textrm{divides}\quad \lambda\binom{n-i}{t-i},\]
apart from a finite number of exceptional $n$ given fixed $k,t,\lambda$.
\end{thm}
Now we are ready to give proofs for our main theorems.

\noindent{Proof of Theorem \ref{final}.} Let $k=2\alpha$. Consider two hypergraphs:

$H_1$ consists of a set of $\alpha$ vertices $K$ such that each $\alpha$-subset of $V(H_1)\setminus K$ forms an edge with $K$.

$H_2$ is a $\frac{k}{2}-(n,k,1)$-design

By the assumption of our theorem and by Theorem \ref{Keevash} a hypergraph $H_2$ exists. Note that 
\[
|E(H_1)|+|E(H_2)|= \binom{n-\alpha}{\alpha}+\frac{\binom{n}{\alpha}}{\binom{2\alpha}{\alpha}}
\]
since each of $\binom{n}{\alpha}$ subsets of vertices in $H_2$ forms exactly one edge, and each such edge is counted $\binom{2\alpha}{\alpha}$ many times.

Observe now that for any bijection $f\colon V(H_1)\to V(H_2)$ the set $K$ is mapped into some set of $\frac{k}{2}$ vertices $K'$ in $H_2$ for wich there exists another set of $\frac{k}{2}$ vertices $K''$ such that $K'\cup K''\in E(H_2)$. Since there is $U\subset V(H_1)\setminus K$ such that $K''=f(U)$ and $U\cup K\in E(H_1)$ for each such $U$, there is no packing of $H_1$ and $H_2$.

The case of odd $k$ is more complicated since we cannot split each edge into two equal pieces. 
Let $t=\lfloor n^{(k-2)/(2k-3)} \rfloor$ Consider now two hypergraphs:

$H_1'$ consists of a complete \hyp on $(k-2)t+1$ vertices and $n-((k-2)t+1)$ independent vertices such that each vertex from the latter part forms an edge with each $(k-1)$-element subset of vertices of the complete part.

$H_2'$ consists of $t$ disjoint copies of $H$, where $H$ is a $(k-1)-(n/t,k,1)$-design.

As before $H_2'$ exists. Hence 
\begin{eqnarray*}
|E(H_1')|+|E(H_2')|&=&\binom{(k-2)t+1}{k} + \binom{(k-2)t+1}{k-1}\Big(n-(k-2)t-1\Big) \\
&&+ \binom{n/t}{k-1}(k-2)t\\
&\leq&c_1t^k+c_2nt^{k-1}-c_3t^k+c_4(\frac{n}{t})^{k-1}t\\
&\leq&cn^{(k^2-k-1)/(2k-3)}
\end{eqnarray*}
for some constant $c,c_1,c_2,c_3,c_4$ depending only on $k$.

Observe that for any bijection $f\colon V(H_1')\to V(H_2')$ at least $k-1$ vertices $\{v_{i_1},\ldots,v_{i_{k-1}}\}$ from the clique in $H_1'$ must be placed onto one of the copies, say $H'$, of $H$ in $H_2'$. Therefore there exists a vertex $x$ in $H'$ such that $\{f(v_{i_1}),\ldots,f(v_{i_{k-1}}), x\}$ is an edge in $H_2'$. Since for every vertex $v\in V(H_1')$ different from $v_{i_1},\ldots,v_{i_{k-1}}$ there is an edge $\{v_{i_1},\ldots,v_{i_{k-1}},v\}$, the case when $v=f^{-1}(x)$ shows that there is no packing of $H_1'$ and $H_2'$.\cbdo
\\

Note that if divisibility conditions from the Theorem \ref{final} are not satisfied we can (in the case of even $k$, the latter is analogous) add $r=r(k)$ isolated vertices (i.e. vertices not belonging to any edge) to $H_2$ and take $H_1$ as a hypergraph that consists of $\lceil r/\alpha\rceil+1$ disjoint $\alpha$-sets $K_1,\ldots,K_{\alpha}$ of vertices and for each $i=1,\ldots,\alpha$ join $K_i$ with each $\alpha$-subset of $V(H_1)\setminus \bigcup_{i=1}^{\alpha}K_i$ and proceed analogously as in the proof. It changes the bound in the theorem, but not the order of the magnitude of this bound.
\\

\noindent {Proof of Theorem \ref{beta}.} Suppose that there exists $0<\beta<k$ such that \ref{beta} is satisfied, but there is no packing of $H_1$ and $H_2$. Let $f$ be any bijection $V(H_1)\to V(H_2)$ and $C$ be any conflict, i.e. a set of vertices of $H_2$ forming an edge and such that their preimages form an edge in $H_1$. Consider any $\beta$-element subset of $C$, say $\{u_1',\ldots,u_{\beta}'\}\subset V(H_2)$ and its correspondng set of preimages $\{u_1,\ldots,u_{\beta}\}\subset V(H_1)$. We want to make a switch of preimages  and some $\{v_1,\ldots,v_{\beta}\}\subset V(H_1)$ by switching the image of $u_i$ with an image $v_i'$ of $v_i$ (in an arbitrary way), such that the new bijection gives fewer conflicts. Obviously $\{v_1,\ldots,v_{\beta}\}\not\subset f^{-1}(C)\}$. The switch is not good if there exist $C_1$ such that 
\[ \{u_1,\ldots,u_{\beta}\}\cup C_1\in E(H_1) \textrm{ and }  \{v_1',\ldots,v_{\beta}'\}\cup f(C_1)\in E(H_2),\]
or there exist $C_2'$ such that 
\[\{u_1',\ldots,u_{\beta}'\}\cup C_2'\in E(H_2) \textrm{ and }  \{v_1,\ldots,v_{\beta}\}\cup f^{-1}(C_2')\in E(H_1).\]
In the first case there can be at most $\Delta_{\beta}(H_1)$ such $C_1$'s  and $f(C_1)$ can be contained in at most $\Delta_{k-\beta}(H_2)$ edges in $H_2$. Since we counted twice $\{u_1,\ldots,u_{\beta}\}= \{v_1,\ldots,v_{\beta}\}$ this forbids
\[\Delta_{\beta}(H_1)\Delta_{k-\beta}(H_2)-1\]
$\beta$-element subsets for switch. Analogously  there can be at most $\Delta_{\beta}(H_2)$ such $C_2'$'s  and $f^{-1}(C_2')$ can be contained in at most $\Delta_{k-\beta}(H_1)$ edges in $H_1$, which forbids
\[\Delta_{k-\beta}(H_1)\Delta_{\beta}(H_2)-1\]
$\beta$-element subsets for switch. Therefore as long as 
\[\Delta_{\beta}(H_1)\Delta_{k-\beta}(H_2)-1+\Delta_{k-\beta}(H_1)\Delta_{\beta}(H_2)-1<\binom{n}{\beta}-\binom{k}{\beta}\]
and there is still a conflict, there is always some $\beta$-element subset among the remaining $\binom{n}{\beta}-\binom{k}{\beta}$ to switch with $\{u_1,\ldots,u_{\al}\}$  yielding a new bijection with fewer conflicts, thus proving the theorem.\cbdo


\begin{thebibliography}{}

\bibitem{A}
N.Alon, Packing of partial designs, Graphs and Combinatorics (1994) 10: 11

\bibitem{Bo}
B. Bollob\'as: Extremal Graph Theory, Academic Press, London-New York, (1978).

\bibitem{BE}
B.~Bollob\'as and S.~E. Eldridge,  \newblock Packing of graphs and applications to
computational complexity,  \newblock{\em J. Combin. Theory Ser. B} 25:105--124, 1978.

\bibitem{K}
P. Keevash. The existence of designs. arXiv:1401.3665, 2014

\bibitem{KSH}
Kostochka, A., Stocker, C. and Hamburger, P. (2013), A Hypergraph Version of a Graph Packing Theorem by Bollobás and Eldridge. J. Graph Theory, 74: 222–235.

\bibitem{N}
P.Naroski, Packing of nonuniform hypergraphs - product and sum of sizes conditions, Discussiones Math. G. Th. 29(3) (2009), 651-656.

\bibitem{PW1}
M. Pil\'sniak, M. Wo\'zniak, A note on packing of two copies of a hypergraph, Discussiones Math. G. Th. 27(1) (2007), 45-49.

\bibitem{PW2}
M. Pil\'sniak, M. Wo\'zniak. On packing of two copies of a hypergraph. Discrete Mathematics and Theoretical Computer Science, DMTCS, 2011, 13 (3), pp.67–74.

\bibitem{RRT}
R\"odl, V., Ruciński, A., Taraz, A. (1999). Hypergraph Packing and Graph Embedding. Combinatorics, Probability and Computing, 8(4), 363-376.

\bibitem{SS}
N.~Sauer and J.~Spencer,  Edge disjoint placement of graphs, 
J. Combin. Theory Ser.~B, 25:295--302, 1978.

\end{thebibliography}
\end{document}